\documentclass[reqno]{amsart}
\usepackage{amsmath,amsthm,amssymb,amscd}

\newtheorem{theorem}{Theorem}[section]

\theoremstyle{definition}
\newtheorem{definition}[theorem]{Definition}

\theoremstyle{remark}

\newtheorem{conj}[theorem]{Conjecture}
\numberwithin{equation}{section}
\allowdisplaybreaks
\begin{document}

%
%
%
%
%
%
%
%
%

\title[Proofs of some conjectures of Keith and Zanello]
 {Proofs of some conjectures of Keith and Zanello on $t$-regular partition}

\author{Ajit Singh}
\address{Department of Mathematics, Indian Institute of Technology Guwahati, Assam, India, PIN- 781039}
\email{ajit18@iitg.ac.in}

\author{Rupam Barman}
\address{Department of Mathematics, Indian Institute of Technology Guwahati, Assam, India, PIN- 781039}
\email{rupam@iitg.ac.in}

\date{January 18, 2022}


\subjclass{Primary 05A17, 11P83, 11F11}

\keywords{$t$-regular partitions; Eta-quotients; modular forms}

\dedicatory{}

\begin{abstract}  For a positive integer $t$, let $b_{t}(n)$ denote the number of $t$-regular partitions of a nonnegative integer $n$. In a recent paper, Keith and Zanello established infinite families of congruences and self-similarity results modulo $2$ for $b_{t}(n)$ for certain values of $t$. Further, they proposed some conjectures on self-similarities of $b_t(n)$ modulo $2$ for certain values of $t$. In this paper, we prove their conjectures on $b_3(n)$ and $b_{25}(n)$. We also prove a self-similarity result for $b_{21}(n)$ modulo $2$ .  
\end{abstract}
\maketitle
\section{Introduction and statement of results} 
 A partition of a positive integer $n$ is any non-increasing sequence of positive integers whose sum is $n$. The number of such partitions of $n$ is denoted by $p(n)$.  
 Let $t$ be a fixed positive integer. A $t$-regular partition of a positive integer $n$ is a partition of $n$ such that none of its part is divisible by $t$. Let $b_{t}(n)$ denote the number of $t$-regular partitions of $n$. The generating function of $b_{t}(n)$ is given by 
\begin{align}\label{gen_fun}
\sum_{n=0}^{\infty}b_{t}(n)q^n=\frac{f_{t}}{f_1},
\end{align}
where $f_k:=(q^k; q^k)_{\infty}=\prod_{j=1}^{\infty}(1-q^{jk})$ and $k$ is a positive integer.
 \par 
 In a very recent paper \cite{Keith2021}, Keith and Zanello studied $t$-regular partition for certain values of $t$. They proved various congruences for $b_t(n)$ modulo $2$ for certain values of $t\leq 28$, and posed several open questions. One of the congruences they proved for $b_3(n)$ is the following:
 \begin{align}\label{new-eqn-20}
 \sum_{n=0}^{\infty}b_{3}(26n+14)q^n\equiv \sum_{n=0}^{\infty}b_{3}(2n)q^{13n}\pmod 2.
 \end{align}
 More generally, they conjectured that:
\begin{conj}\cite[Conjecture 6]{Keith2021}\label{conj1}
For any prime $p>3$, let $\alpha\equiv-24^{-1}\pmod{p^2}$, $0<\alpha<p^2$. It holds for a positive proportion of primes $p$ that 
\begin{align}\label{new-eq-2}
\sum_{n=0}^{\infty}b_{3}(2(pn+\alpha))q^n\equiv \sum_{n=0}^{\infty}b_{3}(2n)q^{pn}\pmod 2.
\end{align}
\end{conj}
The congruence \eqref{new-eqn-20} is a specific case of \eqref{new-eq-2} corresponding to $p=13$. In \cite{singh-barman}, we proved a specific case of \eqref{new-eq-2} corresponding to $p=17$. The aim of this article is to prove two conjectures of Keith and Zanello. Our first theorem confirms Conjecture \ref{conj1}.
\begin{theorem} \label{thm1}
Conjecture \ref{conj1} is true.
\end{theorem}
Keith and Zanello also studied $2$-divisibility of $b_{25}(n)$ and proved several congruences for primes $p\equiv 11, 13,17,19 \pmod{20}$ and $p\equiv 31,39 \pmod{40}$. 
To be specific, if $p\equiv 11, 13,17,19 \pmod{20}$ is prime, then they proved that
\begin{align*}
b_{25}(8(p^2n+kp-3\cdot 4^{-1})+5)\equiv 0\pmod{2}
\end{align*}
for all $1\leq k<p$, where $3\cdot 4^{-1}$ is taken modulo $p^2$.  
Further, they conjectured the following:
\begin{conj}\cite[Conjecture 28]{Keith2021}\label{conj3}
	For a positive proportion of primes $p$, it holds that
	\begin{align*}
	\sum_{n=0}^{\infty} b_{25}(2 p n+\alpha) q^{n} \equiv q^{\beta} \sum_{n=0}^{\infty} b_{25}(2 n+1) q^{p n}\pmod2,
	\end{align*}
	for some $\alpha$ and $\beta$ depending on $p$.
\end{conj}
Our second theorem confirms Conjecture \ref{conj3}.
\begin{theorem}\label{thm3}
	Conjecture \ref{conj3} is true.
\end{theorem}
Next, we prove a self-similarity result for $b_{21}(n)$ modulo 2. More precisely, we prove the following theorem:
\begin{theorem}\label{thm4}
	For a positive proportion of primes $p$, it holds that
	\begin{align*}
	\sum_{n=0}^{\infty}b_{21}(pn+11\gamma+1)q^{n}\equiv\sum_{n=0}^{\infty}b_{21}(n+1)q^{pn}\pmod 2,
	\end{align*}
	for some $\gamma$ depending on $p$.
\end{theorem}
\section{Preliminaries}
We recall some definitions and basic facts on modular forms. For more details, see for example \cite{koblitz1993, ono2004}. We first define the matrix groups 
\begin{align*}
\text{SL}_2(\mathbb{Z}) & :=\left\{\begin{bmatrix}
a  &  b \\
c  &  d      
\end{bmatrix}: a, b, c, d \in \mathbb{Z}, ad-bc=1
\right\},\\
\Gamma_{0}(N) & :=\left\{
\begin{bmatrix}
a  &  b \\
c  &  d      
\end{bmatrix} \in \text{SL}_2(\mathbb{Z}) : c\equiv 0\pmod N \right\},
\end{align*}
\begin{align*}
\Gamma_{1}(N) & :=\left\{
\begin{bmatrix}
a  &  b \\
c  &  d      
\end{bmatrix} \in \Gamma_0(N) : a\equiv d\equiv 1\pmod N \right\},
\end{align*}
and 
\begin{align*}\Gamma(N) & :=\left\{
\begin{bmatrix}
a  &  b \\
c  &  d      
\end{bmatrix} \in \text{SL}_2(\mathbb{Z}) : a\equiv d\equiv 1\pmod N, ~\text{and}~ b\equiv c\equiv 0\pmod N\right\},
\end{align*}
where $N$ is a positive integer. A subgroup $\Gamma$ of $\text{SL}_2(\mathbb{Z})$ is called a congruence subgroup if $\Gamma(N)\subseteq \Gamma$ for some $N$. The smallest $N$ such that $\Gamma(N)\subseteq \Gamma$
is called the level of $\Gamma$. For example, $\Gamma_0(N)$ and $\Gamma_1(N)$
are congruence subgroups of level $N$. 
\par Let $\mathbb{H}:=\{z\in \mathbb{C}: \text{Im}(z)>0\}$ be the upper half of the complex plane. The group $$\text{GL}_2^{+}(\mathbb{R})=\left\{\begin{bmatrix}
a  &  b \\
c  &  d      
\end{bmatrix}: a, b, c, d\in \mathbb{R}~\text{and}~ad-bc>0\right\}$$ acts on $\mathbb{H}$ by $\begin{bmatrix}
a  &  b \\
c  &  d      
\end{bmatrix} z=\displaystyle \frac{az+b}{cz+d}$.  
We identify $\infty$ with $\displaystyle\frac{1}{0}$ and define $\begin{bmatrix}
a  &  b \\
c  &  d      
\end{bmatrix} \displaystyle\frac{r}{s}=\displaystyle \frac{ar+bs}{cr+ds}$, where $\displaystyle\frac{r}{s}\in \mathbb{Q}\cup\{\infty\}$.
This gives an action of $\text{GL}_2^{+}(\mathbb{R})$ on the extended upper half-plane $\mathbb{H}^{\ast}=\mathbb{H}\cup\mathbb{Q}\cup\{\infty\}$. 
Suppose that $\Gamma$ is a congruence subgroup of $\text{SL}_2(\mathbb{Z})$. A cusp of $\Gamma$ is an equivalence class in $\mathbb{P}^1=\mathbb{Q}\cup\{\infty\}$ under the action of $\Gamma$.
\par The group $\text{GL}_2^{+}(\mathbb{R})$ also acts on functions $f: \mathbb{H}\rightarrow \mathbb{C}$. In particular, suppose that $\gamma=\begin{bmatrix}
a  &  b \\
c  &  d      
\end{bmatrix}\in \text{GL}_2^{+}(\mathbb{R})$. If $f(z)$ is a meromorphic function on $\mathbb{H}$ and $\ell$ is an integer, then define the slash operator $|_{\ell}$ by 
$$(f|_{\ell}\gamma)(z):=(\text{det}~{\gamma})^{\ell/2}(cz+d)^{-\ell}f(\gamma z).$$
\begin{definition}
	Let $\Gamma$ be a congruence subgroup of level $N$. A holomorphic function $f: \mathbb{H}\rightarrow \mathbb{C}$ is called a modular form with integer weight $\ell$ on $\Gamma$ if the following hold:
	\begin{enumerate}
		\item We have $$f\left(\displaystyle \frac{az+b}{cz+d}\right)=(cz+d)^{\ell}f(z)$$ for all $z\in \mathbb{H}$ and all $\begin{bmatrix}
		a  &  b \\
		c  &  d      
		\end{bmatrix} \in \Gamma$.
		\item If $\gamma\in \text{SL}_2(\mathbb{Z})$, then $(f|_{\ell}\gamma)(z)$ has a Fourier expansion of the form $$(f|_{\ell}\gamma)(z)=\displaystyle\sum_{n\geq 0}a_{\gamma}(n)q_N^n,$$
		where $q_N:=e^{2\pi iz/N}$.
	\end{enumerate}
In addition, if $a_{\gamma}(0)=0$ for all $\gamma \in \text{SL}_2(\mathbb{Z})$, then $f$ is called a cusp form.
\end{definition}
For a positive integer $\ell$, the complex vector space of modular forms (resp. cusp forms) of weight $\ell$ with respect to a congruence subgroup $\Gamma$ is denoted by $M_{\ell}(\Gamma)$ (resp. $S_{\ell}(\Gamma)$).
\begin{definition}\cite[Definition 1.15]{ono2004}
	If $\chi$ is a Dirichlet character modulo $N$, then we say that a modular form $f\in M_{\ell}(\Gamma_1(N))$ (resp. $S_{\ell}(\Gamma_1(N))$) has Nebentypus character $\chi$ if
	$$f\left( \frac{az+b}{cz+d}\right)=\chi(d)(cz+d)^{\ell}f(z)$$ for all $z\in \mathbb{H}$ and all $\begin{bmatrix}
	a  &  b \\
	c  &  d      
	\end{bmatrix} \in \Gamma_0(N)$. The space of such modular forms (resp. cusp forms) is denoted by $M_{\ell}(\Gamma_0(N), \chi)$ (resp. $S_{\ell}(\Gamma_0(N), \chi)$). 
\end{definition}
In this paper, the relevant modular forms are those that arise from eta-quotients. Recall that the Dedekind eta-function $\eta(z)$ is defined by
\begin{align*}
\eta(z):=q^{1/24}(q;q)_{\infty}=q^{1/24}\prod_{n=1}^{\infty}(1-q^n),
\end{align*}
where $q:=e^{2\pi iz}$ and $z\in \mathbb{H}$. A function $f(z)$ is called an eta-quotient if it is of the form
\begin{align*}
f(z)=\prod_{\delta\mid N}\eta(\delta z)^{r_\delta},
\end{align*}
where $N$ is a positive integer and $r_{\delta}$ is an integer. We now recall two theorems from \cite[p. 18]{ono2004} which are very useful in checking modularity of eta-quotients.
\begin{theorem}\cite[Theorem 1.64]{ono2004}\label{thm_ono1} If $f(z)=\prod_{\delta\mid N}\eta(\delta z)^{r_\delta}$ 
	is an eta-quotient such that $\ell=\frac{1}{2}\sum_{\delta\mid N}r_{\delta}\in \mathbb{Z}$, 
	$$\sum_{\delta\mid N} \delta r_{\delta}\equiv 0 \pmod{24}$$ and
	$$\sum_{\delta\mid N} \frac{N}{\delta}r_{\delta}\equiv 0 \pmod{24},$$
	then $f(z)$ satisfies $$f\left( \frac{az+b}{cz+d}\right)=\chi(d)(cz+d)^{\ell}f(z)$$
	for every  $\begin{bmatrix}
	a  &  b \\
	c  &  d      
	\end{bmatrix} \in \Gamma_0(N)$. Here the character $\chi$ is defined by $\chi(d):=\left(\frac{(-1)^{\ell} s}{d}\right)$, where $s:= \prod_{\delta\mid N}\delta^{r_{\delta}}$. 
\end{theorem}
Suppose that $f$ is an eta-quotient satisfying the conditions of Theorem \ref{thm_ono1} and that the associated weight $\ell$ is a positive integer. If $f(z)$ is holomorphic (resp. vanishes) at all of the cusps of $\Gamma_0(N)$, then $f(z)\in M_{\ell}(\Gamma_0(N), \chi)$ (resp. $S_{\ell}(\Gamma_0(N), \chi)$). The following theorem gives the necessary criterion for determining orders of an eta-quotient at cusps.
\begin{theorem}\cite[Theorem 1.65]{ono2004}\label{thm_ono2}
	Let $c, d$ and $N$ be positive integers with $d\mid N$ and $\gcd(c, d)=1$. If $f$ is an eta-quotient satisfying the conditions of Theorem \ref{thm_ono1} for $N$, then the 
	order of vanishing of $f(z)$ at the cusp $\frac{c}{d}$ 
	is $$\frac{N}{24}\sum_{\delta\mid N}\frac{\gcd(d,\delta)^2r_{\delta}}{\gcd(d,\frac{N}{d})d\delta}.$$
\end{theorem}
\par We next recall the definition of Hecke operators.
Let $m$ be a positive integer and $f(z) = \sum_{n=0}^{\infty} a(n)q^n \in M_{\ell}(\Gamma_0(N),\chi)$. Then the action of Hecke operator $T_m$ on $f(z)$ is defined by 
\begin{align*}
f(z)|T_m := \sum_{n=0}^{\infty} \left(\sum_{d\mid \gcd(n,m)}\chi(d)d^{\ell-1}a\left(\frac{nm}{d^2}\right)\right)q^n.
\end{align*}
In particular, if $m=p$ is prime, we have 
\begin{align*}
f(z)|T_p := \sum_{n=0}^{\infty} \left(a(pn)+\chi(p)p^{\ell-1}a\left(\frac{n}{p}\right)\right)q^n.
\end{align*}
We adopt the convention that $a(n/p)=0$ when $p\nmid n$. 
\par We finally recall a result of Serre \cite{serre1} (also see \cite[Proposition 4.2]{Treneer2006}) about the action of Hecke operator on cusp forms. For a number field $K$, let $\mathcal{O}_{K}$ denote its ring of integers.
\begin{theorem}\cite[Exercise 6.4]{serre1}\label{Serre1}
Suppose that
$$f(z)=\sum_{n=1}^{\infty}a(n)q^n\in S_k(\Gamma_{0}(N),\chi)$$
has coefficients in $\mathcal{O}_{K}$, and $M$ is a positive integer. Furthermore, suppose that $k>1$. Then a positive proportion of the primes $p\equiv-1\pmod{MN}$ have the property that
$$f(z) \mid T_{p} \equiv 0 \pmod M.$$
\end{theorem}
\section{Proof of Theorem \ref{thm1}}
\begin{proof}
		We first recall the following even-odd disection of the $3$-regular partitions \cite[(6)]{Keith2021}:
	\begin{align*}
	\sum_{n=0}^{\infty}b_3(n)q^n=\frac{f_3}{f_1}\equiv\frac{f_1^8}{f_3^2}+q\frac{f_3^{10}}{f_1^4}\pmod 2.
	\end{align*}
	Extracting the terms with even powers of $q$, we obtain
	\begin{align}\label{eqn-new-3}
	\sum_{n=0}^{\infty}b_3(2n)q^n\equiv \frac{f_1^4}{f_3}\pmod 2.
	\end{align} 
	Let 
	\begin{align*}
	A(z) := \prod_{n=1}^{\infty} \frac{(1-q^{24n})^2}{(1-q^{48n})} = \frac{\eta^2(24z)}{\eta(48z)}. 
	\end{align*}
	Then using the binomial theorem we have 
	\begin{align*}
	A(z) = \frac{\eta^2(24z)}{\eta(48z)} \equiv 1 \pmod {2}.
	\end{align*}
	Define $B(z)$ by
	\begin{align*}
	B(z):= \left(\frac{\eta^{4}(24z)}{\eta(72z)}\right)A(z)
	=\frac{\eta^{6}(24z)}{\eta(72z)\eta(48z)}.
	\end{align*}
	Modulo $2$, we have
	\begin{align}\label{thm1.1}
	B(z)\equiv\frac{\eta^{4}(24z)}{\eta(72z)} =q\frac{(q^{24}; q^{24})_{\infty}^4}{(q^{72}; q^{72})_{\infty}}.
	\end{align}
	Combining \eqref{eqn-new-3} and \eqref{thm1.1}, we obtain 
	\begin{align}\label{thm1.2}
	B(z) \equiv \sum_{n=0}^{\infty}b_3(2n)q^{24n+1} \pmod {2}.
	\end{align}
	Now, $B(z)$ is an eta-quotient with $N =3456$. We next prove that $B(z)$ is a modular form. We know that the cusps of $\Gamma_{0}(3456)$ are represented by fractions $\frac{c}{d}$, where $d\mid 3456$ and $\gcd(c, d)=1$. By Theorem \ref{thm_ono2}, we find that $B(z)$ vanishes at a cusp $\frac{c}{d}$ if and only if
	\begin{align*}
	L:=12\frac{\gcd(d,24)^2}{\gcd(d,48)^2}-\frac{2}{3}\frac{\gcd(d,72)^2}{\gcd(d,48)^2}-1> 0.
	\end{align*}
	We now consider the following four cases according to the divisors of $3456$ and find the values of $G_1:=\frac{\gcd(d,24)^2}{\gcd(d,48)^2}$ and $G_2:=\frac{\gcd(d,72)^2}{\gcd(d,48)^2}$. Let $d$ be a divisor of $N=3456$. \\
	Case (i). For $d=2^{r_1}3^{r_2}$, where $0\leq r_1\leq 3$ and $0\leq r_2\leq 1$, we find that $G_1=G_2=1$. Hence, $L>0$.\\
	Case (ii). For $d=2^{r_1}3^{r_2}$, where $0\leq r_1\leq 3$ and $2\leq r_2\leq 3$, we find that $G_1=1$ and $G_2=9$. Hence, $L>0$.\\
	Case (iii). For $d=2^{r_1}3^{r_2}$, where $4\leq r_1\leq 7$ and $0\leq r_2\leq 1$, we find that $G_1=G_2=1/4$. Hence, $L>0$.\\
	Case (iv). For $d=2^{r_1}3^{r_2}$, where $4\leq r_1\leq 7$ and $2\leq r_2\leq 3$, we find that $G_1=1/4$ and $G_2=9/4$. Hence, $L>0$.
	\par
	Thus, $B(z)$ vanishes at every cusp $\frac{c}{d}$. Using Theorem \ref{thm_ono1}, we find that the weight of $B(z)$ is equal to $2$. Also, the associated character for $B(z)$ is given by $\chi_1=(\frac{2^{11} 3^{3}}{\bullet})$.
	This proves that $B(z) \in S_{2}(\Gamma_{0}(3456), \chi_1)$. Also, the Fourier coefficients of $B(z)$ are all integers. Hence by Theorem \ref{Serre1}, a positive proportion of the primes $p\equiv-1\pmod{6912}$ have the property that 
	\begin{align}\label{thm1.3}
	B(z) \mid T_{p} \equiv 0 \pmod 2.
	\end{align}	
	Let $B(z)=\sum_{n=1}^{\infty}a(n)q^n$. Then, \eqref{thm1.2} yields
	\begin{align}\label{thm1.5}
	\sum_{n=1}^{\infty}b_3\left(\frac{2(n-1)}{24}\right)q^n\equiv\sum_{n=1}^{\infty}a(n)q^n\pmod 2.
	\end{align}
	Now, from \eqref{thm1.3} we obtain
	\begin{align*}
	B(z) \mid T_{p}=\sum_{n=1}^{\infty}(a(pn)+p\chi_1(p) a(n/p))q^n\equiv 0\pmod 2
	\end{align*}	
	which yields 
	\begin{align}\label{thm1.4}
	\sum_{n=1}^{\infty}a(pn)q^n \equiv  \sum_{n=1}^{\infty}a(n/p)q^n\pmod 2.
	\end{align}
Combining \eqref{thm1.5} and \eqref{thm1.4} we find that 
	\begin{align*}
	\sum_{n=1}^{\infty}b_3\left(\frac{2(pn-1)}{24}\right)q^n&\equiv\sum_{n=1}^{\infty}b_3\left(\frac{2(n/p-1)}{24}\right)q^n\pmod 2\\
	&\equiv\sum_{n=1}^{\infty}b_3\left(\frac{2(n-1)}{24}\right)q^{pn}\pmod 2\\
	&\equiv\sum_{n=0}^{\infty}b_3\left(\frac{2n}{24}\right)q^{pn+p}\pmod 2.
	\end{align*}
	Multiplying both sides by $q^{-p}$ we obtain
	\begin{align*}
	&\sum_{n=p}^{\infty}b_3\left(\frac{2(pn-1)}{24}\right)q^{n-p}\equiv\sum_{n=0}^{\infty}b_3\left(\frac{2n}{24}\right)q^{pn}\pmod 2
	\end{align*}
	which yields 
	\begin{align*}
	&\sum_{n=0}^{\infty}b_3\left(\frac{2(pn+p^2-1)}{24}\right)q^n\equiv\sum_{n=0}^{\infty}b_3\left(\frac{2n}{24}\right)q^{pn}\pmod 2.
	\end{align*}
	Let $\alpha=\frac{p^2-1}{24}$. Since, $p\equiv-1\pmod{6912}$, so $\alpha$ is a positive integer, and $\alpha\equiv -24^{-1}\pmod {p^2}$, $0<\alpha<p^2$.
	Replacing  $n$ by $24n$ and then substituting $q^{24}$ by $q$ we get
	\begin{align*}
	\sum_{n=0}^{\infty}b_3(2(pn+\alpha))q^{n}\equiv\sum_{n=0}^{\infty}b_3(2n)q^{pn}\pmod 2.
	\end{align*}
	This completes the proof of the theorem.
\end{proof}
\section{Proof of Theorem \ref{thm3}}
\begin{proof} 
Putting $t=25$ in \eqref{gen_fun} we have
\begin{align}\label{thm3.01}
\sum_{n=0}^{\infty}b_{25}(n)q^n=\frac{f_{25}}{f_1}.
\end{align}
We use identity \cite[(4)]{Judge2018}, namely
\begin{align*}
f_1f_5\equiv f_1^6+qf_5^6 \pmod 2.
\end{align*}
Dividing both sides by $f_1^2$ we obtain
\begin{align}\label{thm3.02}
\frac{f_5}{f_1}\equiv f_1^4+q\frac{f_5^6}{f_1^2}\pmod 2.
\end{align}
Therefore, by \eqref{thm3.01} and \eqref{thm3.02} we have
\begin{align*}
\sum_{n=0}^{\infty}b_{25}(n)q^n&=\frac{f_{25}}{f_1}=\frac{f_{25}}{f_5}\frac{f_{5}}{f_1}\\
&\equiv f_1^4f_5^4+q^6f_5^4\frac{f_{25}^6}{f_1^2}+q\frac{f_5^{10}}{f_1^2}+q^5f_1^4\frac{f_{25}^6}{f_5^2}\pmod 2.
\end{align*}
Extracting the terms involving $q^{2n+1}$, and then dividing by $q$ and replacing $q^2$ by $q$, we find that 
\begin{align*}
\sum_{n=0}^{\infty}b_{25}(2n+1)q^n&\equiv \frac{f_5^5 }{f_1}+q^2 \frac{f_1^2 f_{25}^3}{f_5}\pmod2\\
&\equiv f_{1}^{4} f_{5}^{4}+q \frac{f_{5}^{10}}{f_{1}^{2}}+q^{2} f_{1}^{2} f_{5}^{4} f_{25}^{2}+q^{7} \frac{f_{1}^{2} f_{25}^{8}}{f_{5}^{2}}\pmod 2. 
\end{align*}
Extracting the terms involving $q^{2n}$, we obtain
\begin{align}\label{thm3.1}
\sum_{n=0}^{\infty}b_{25}(4n+1)q^{2n}\equiv f_{2}^{2} f_{10}^{2}+q^2 f_{2} f_{10}^{2} f_{50}\pmod 2.
\end{align}
Define $F(z)$ by
\begin{align}\label{thm3.2}
F(z)&:=\eta^2(2z)\eta^2(10z)+\eta(2z)\eta^2(10z)\eta(50z).
\end{align}
Combining \eqref{thm3.1} and \eqref{thm3.2}, we obtain 
\begin{align}\label{thm3.4}
F(z) \equiv \sum_{n=0}^{\infty}b_{25}(4n+1)q^{2n+1} \pmod {2}.
\end{align}
Now using Theorems \ref{thm_ono1} and \ref{thm_ono2}, we find that $\eta^2(2z)\eta^2(10z)\in S_{2}\left(\Gamma_{0}(100), \chi_3\right)$ and $\eta(2z)\eta^2(10z)\eta(50z)\in S_{2}\left(\Gamma_{0}(100), \chi_3\right)$ for some Nebentypus character $\chi_3$ and hence $F(z)\in S_{2}\left(\Gamma_{0}(100), \chi_3\right)$. Also, the Fourier coefficients of $F(z)$ are all integers. Hence by Theorem \ref{Serre1}, a positive proportion of the primes $p\equiv-1\pmod{200}$ have the property that 
\begin{align}\label{thm3.5}
F(z) \mid T_{p} \equiv 0 \pmod 2.
\end{align}	
Let $F(z)=\sum_{n=1}^{\infty}d(n)q^n$. Then, \eqref{thm3.4} yields
\begin{align}\label{thm3.7}
\sum_{n=1}^{\infty}b_{25}(2(n-1)+1)q^n\equiv\sum_{n=1}^{\infty}d(n)q^n\pmod 2.
\end{align}
Now, from \eqref{thm3.5} we obtain
\begin{align*}
F(z) \mid T_{p}=\sum_{n=1}^{\infty}(d(pn)+p\chi_3(p) d(n/p))q^n\equiv 0\pmod 2
\end{align*}	
which yields 
\begin{align}\label{thm3.6}
\sum_{n=1}^{\infty}d(pn)q^n\equiv \sum_{n=1}^{\infty}d(n/p)q^n\pmod 2.
\end{align}
Combining \eqref{thm3.7} and \eqref{thm3.6} we find that  
\begin{align*}
\sum_{n=1}^{\infty}b_{25}(2pn-1)q^n&\equiv\sum_{n=1}^{\infty}b_{25}(2(n/p-1)+1)q^n\pmod 2\\
&\equiv\sum_{n=1}^{\infty}b_{25}(2(n-1)+1)q^{pn}\pmod 2\\
&\equiv\sum_{n=0}^{\infty}b_{25}(2n+1)q^{pn+p}\pmod 2.
\end{align*}
Replacing $n$ by $n+1$ on the left side and then dividing both sides by $q$ we obtain
\begin{align*}
\sum_{n=0}^{\infty}b_{25}(2pn+\alpha)q^{n}\equiv q^{\beta}\sum_{n=0}^{\infty}b_{25}(2n+1)q^{pn}\pmod 2,
\end{align*}
where $\alpha=2p-1$ and $\beta=p-1$. This completes the proof of the theorem.
\end{proof}
\section{Proof of Theorem \ref{thm4}}
\begin{proof} 
	We begin with the identity \cite[Section 7]{Keith2021}, namely
	\begin{align}\label{thm4.1}
	\sum_{n=0}^{\infty}b_{21}(4n+1)q^n\equiv\frac{f_{3}^4}{f_1}\pmod 2.
	\end{align}
	Let 
	\begin{align*}
	G(z) := \prod_{n=1}^{\infty} \frac{(1-q^{24n})^2}{(1-q^{48n})} = \frac{\eta^2(24z)}{\eta(48z)}. 
	\end{align*}
	Then using the binomial theorem we have 
	\begin{align*}
	G(z) = \frac{\eta^2(24z)}{\eta(48z)} \equiv 1 \pmod {2}.
	\end{align*}
	Define $H(z)$ by
	\begin{align*}
	H(z):= \left(\frac{\eta^{4}(72z)}{\eta(24z)}\right)G(z)
	=\frac{\eta^{4}(72z)\eta(24z)}{\eta(48z)}.
	\end{align*}
	Modulo $2$, we have
	\begin{align}\label{thm4.3}
	H(z)\equiv\frac{\eta^{4}(72z)}{\eta(24z)} =q^{11}\frac{(q^{72}; q^{72})_{\infty}^4}{(q^{24}; q^{24})_{\infty}}.
	\end{align}
	Combining \eqref{thm4.1} and \eqref{thm4.3}, we obtain 
	\begin{align}\label{thm4.4}
	H(z) \equiv \sum_{n=0}^{\infty}b_{21}(4n+1)q^{24n+11} \pmod {2}.
	\end{align}
	Now using Theorems \ref{thm_ono1} and \ref{thm_ono2}, we find that $H(z)\in S_{2}\left(\Gamma_{0}(3456), \chi_2\right)$ for some Nebentypus character $\chi_2$. Also, the Fourier coefficients of $H(z)$ are all integers. Hence by Theorem \ref{Serre1}, a positive proportion of the primes $p\equiv-1\pmod{6912}$ have the property that 
	\begin{align}\label{thm4.5}
	H(z) \mid T_{p} \equiv 0 \pmod 2.
	\end{align}	
	Let $H(z):=\sum_{n=1}^{\infty}c(n)q^n$. Then, \eqref{thm4.4} yields
	\begin{align}\label{thm4.7}
	\sum_{n=1}^{\infty}b_{21}\left(\frac{4(n-11)}{24}+1\right)q^n\equiv\sum_{n=1}^{\infty}c(n)q^n\pmod 2.
	\end{align}
	Now, from \eqref{thm4.5} we obtain
	\begin{align*}
	H(z) \mid T_{p}=\sum_{n=1}^{\infty}(c(pn)+p\chi_2(p) c(n/p))q^n\equiv 0\pmod 2
	\end{align*}	
	which yields 
	\begin{align}\label{thm4.6}
	\sum_{n=1}^{\infty}c(pn)q^n\equiv  \sum_{n=1}^{\infty}c(n/p)q^n \pmod 2.
	\end{align}
	Combining \eqref{thm4.7} and \eqref{thm4.6} we find that  
	\begin{align*}
	\sum_{n=1}^{\infty}b_{21}\left(\frac{4(pn-11)}{24}+1\right)q^n&\equiv\sum_{n=1}^{\infty}b_{21}\left(\frac{4(n/p-11)}{24}+1\right)q^n\pmod 2\\
	&\equiv\sum_{n=1}^{\infty}b_{21}\left(\frac{4(n-11)}{24}+1\right)q^{pn}\pmod 2\\
	&\equiv\sum_{n=0}^{\infty}b_{21}\left(\frac{4n}{24}+1\right)q^{pn+11p}\pmod 2.
	\end{align*}
	Multiplying both sides by $q^{-11p}$ we obtain
	\begin{align*}
	&\sum_{n=11p}^{\infty}b_{21}\left(\frac{pn-11}{6}+1\right)q^{n-11p}\equiv\sum_{n=0}^{\infty}b_{21}\left(\frac{n}{6}+1\right)q^{pn}\pmod 2
	\end{align*}
	which yields 
	\begin{align*}
	&\sum_{n=0}^{\infty}b_{21}\left(\frac{pn+11(p^2-1)}{6}+1\right)q^n\equiv\sum_{n=0}^{\infty}b_{21}\left(\frac{n}{6}+1\right)q^{pn}\pmod 2.
	\end{align*}
	Let $\gamma=\frac{p^2-1}{6}$. Since $p\equiv-1\pmod{6912}$, so $\gamma$ is a positive integer, and $\gamma\equiv -6^{-1}\pmod {p^2}$, $0<\gamma<p^2$.
	Replacing  $n$ by $6n$ and then substituting $q^{6}$ by $q$ we get
	\begin{align*}
	\sum_{n=0}^{\infty}b_{21}(pn+11\gamma+1)q^{n}\equiv\sum_{n=0}^{\infty}b_{21}(n+1)q^{pn}\pmod 2.
	\end{align*}
	This completes the proof of the theorem.
\end{proof}
\section{Acknowledgements}
We are extremely grateful to Professor Ken Ono for many helpful discussions while working on this project.

\end{document}